\documentclass{articles_perso}
\title{Instability of the magnetohydrodynamics system at vanishing Reynolds number}
\date{}

\begin{document}
\begin{abstract}
  The aim of this note is to study the dynamo properties of the
  magnetohydrodynamics system at vanishing \(R_m\). Improving the analysis
  in \cite{gerardvaret2006osi}, we shall establish a generic Lyapunov
  instability result.
\end{abstract}
\section{The induction equation}
The incompressible Magnetohydrodynamics (MHD) system is a classical model
for conducting fluids. It is derived from both Navier-Stokes equation for
the fluid movement under the magnetic field and the Maxwell equations with
Ohm's law for the evolution of the magnetic field generated by the fluid
movement. It is in dimensionless form:
\begin{equation}
  \label{inimhd}
  \left\{\begin{array}{l}
      \partial_tu + u\cdot\nabla u + \nabla p - \frac1{R_e}\Delta u = (\nabla\wedge b)\wedge b + f,\\
      \partial_tb-\nabla\wedge(u\wedge b) -\frac1{R_m}\Delta b = 0,\\
      \nabla\cdot u = 0,\\
      \nabla\cdot b = 0,
    \end{array}\right.
\end{equation}
where \(R_e\) and \(R_m\) are the hydrodynamic and magnetic Reynolds
numbers. In this system, \(u = u(t,\theta)\) and \(b = b(t,\theta)\) are
respectively the velocity field of the fluid and the magnetic field,
vector-valued functions in \(\mb R^3\), and \(f\) is an additional forcing
term (anything that is not the electromagnetic force). To lighten
notations, we will assume in this note that \(R_e=1\) and
\(R_m=\varepsilon\), where \(\varepsilon\) is a small parameter. Note that
the Reynolds magnetic number -- that is the one that corresponds to the
definition commonly used -- is defined by \(R_m = \frac{LV}{\nu_m}\), where
\(L\) and \(V\) are respectively the typical variation length and amplitude
of the velocity field \(u\). This means that for a given time-and-space
scale for the velocity field, we consider large magnetic diffusivity (or
vanishing conductivity) of the fluid. The reader may refer to
\cite{fauve2003dynamo} to have a complete explanation on how the above MHD
system is computed. One may say it is the base system for describing MHD in
the case of a conducting viscous fluid.

Our aim here, related to dynamo theory is to deal with the stability of
solutions of this last system, and more precisely with the instability of
solutions of the form
\begin{equation}
  \label{u0}
  (u,b) = (u(t,x),0).
\end{equation}
The physical problem behind the study of the stability of such solutions is
to understand the generation of magnetic field by the fluid. The fluid
``gives'' energy to the magnetic field by the term \(\nabla\wedge(u\wedge
b)\), which should compensate the dissipation term in the equation. We know
for long that this kind of phenomena occurs in planets and stars like the
Sun or the Earth and make them act like a magnet (See
\cite{larmor1919could} where the problem was first posed).

There have been many results in the past century in dynamo theory. First
results where ``negative'' results in the sense that they gave condition
for which fluids cannot generate a dynamo effect, i.e. this corresponds to
situations where solutions of the form \eqref{u0} were actually
stable. Those results spread the idea that the velocity field has to be
chaotic enough to expect a dynamo effect (See M.M.Vishik,
\cite{vishik1989mfg}). The knowledge in the domain concerning stability of
the MHD system has been summed up in \cite{gilbert2003dt}.

Then positive results began to emerge concerning kinetic dynamos, that is a
dynamo where we consider only the induction equation in the system:
physically, this corresponds to the case where the Laplace force is
neglected in the fluid movement, i.e. the magnetic field has no
retro-action on the fluid. This case is much simpler since it corresponds
to only a linear case. More recent results gave positive result with the
full MHD system, which is non-linear.

In this note, we will focus on a particular mechanism for the generation of
magnetic field, that is the ``alpha-effect'', which was first introduced by
E.N. Parker \cite{parker1955hydromagnetic}. This mechanism is based on
scale separation: we formally decompose the fields into two parts, a
fluctuating one, evolving on small lengths, and a mean one evolving on much
bigger ones. The idea is that the small-scale field can have an effect on
the big-scale one through the mean part of the crossed-term
\(\nabla\wedge(u\wedge b)\) in the induction equation and may create
instability on the large-scale. Note that this kind of mechanism has been
experimentally confirmed much recently by R. Stieglitz and U. M\"uller
\cite{stieglitz2001experimental}.

Three regimes of \(R_m\) can be identified in the MHD system:
\begin{description}
\item[The intermediate regime where \(R_m\) is of order \(1\).] We
  considered this regime in a former article \cite{bouya2011instability}.
\item[The regime of large \(R_m\).]  In that case, one can expect a dynamo
  at the scale of the velocity field, since the diffusive effect is very
  weak at such scale. A related question is whether the dynamo is ``fast''
  or ``slow'', that is whether the best growing mode decreases to \(0\) as
  \(\tend{R_m}{+\infty}\) or not. So far this regime has been mainly
  studied formally and numerically (up to \(R_m\sim10^5\)) in physics
  \cite{otani1993fast,soward1987fast,soward1994role,alexakis2011searching},
  in particular in the case of ABC-like flows, and the result seems to be
  difficult to predict in the considered cases. This regime will be the
  object of a forthcoming paper.
\item[The regime of small \(R_m\).] This regime is relevant to laboratory
  experiments, and is the focus of the present note. More precisely, we use
  the \(\alpha\)-effect mechanism to obtain a linearized instability result
  for system \eqref{inimhd}, around steady states \((u,b) = (U(\theta),
  0)\). Then, we conclude to a Lyapunov instability result for the full
  non-linear system.
\end{description}

The base flows we consider satisfy
\begin{equation}
  U\in\mc H^\infty(\mb T^3)^3
\end{equation}
and
\begin{equation}
  \int_{\mb T^3}U=0,\text{ and}\qquad \nabla_\theta\cdot U=0.
\end{equation}
We denote by \(\ms P\) the set of fields \(U\) of this form. \(\ms P\) is a
Fr\'echet space for the associated semi-norms
\begin{equation}
  \|U\|_m^2=\sum_{\xi\in\mb Z^4}|\xi|^{2m}|\hat U(\xi)|^2.
\end{equation}
We have
\begin{equation}
  u = U + u', \quad 
  b = 0 + b',
\end{equation}
focusing on the growth of the perturbations \(u',b'\). As we shall see, it is
also convenient to do a rescaling in time: \(t'= \varepsilon^3t\). Dropping
the primes, we end up with:
\begin{equation}
  \label{fullmhd}
  \left\{\begin{array}{l}
      \partial_tb - \varepsilon^{-3}\nabla\wedge(U\wedge b) - \varepsilon^{-3}\nabla\wedge(u\wedge b) - \varepsilon^{-4}\Delta b = 0,\\
      \varepsilon^3\partial_tu + u\cdot\nabla U +U\cdot\nabla u +u\cdot\nabla u + \nabla p - \Delta u = (\nabla\wedge b)\wedge b\\
      \nabla\cdot b= 0\\
      \nabla\cdot u = 0
    \end{array}\right.
\end{equation}
As mentioned before, a key point is the understanding of the linearized
induction equation,
\begin{equation}
  \label{mhd}
  \partial_tb - \varepsilon^{-3}\nabla\wedge(U\wedge b) - \varepsilon^{-4}\Delta b = 0,
\end{equation}
Motivated by the analysis of the \(\alpha\)-effect, we will look for an
exponentially growing solution of that equation based on scale separation:
\begin{equation}
  b^\varepsilon(t,\theta) = e^{\lambda^\varepsilon t}(\bar b(\varepsilon^2\theta) + \varepsilon\tilde b(\varepsilon^2\theta,\theta)),
\end{equation}
where \(\bar b = \bar b(x)\) and \(\tilde b = \tilde b(x,\theta)\) depend on the original
small-scale variable \(\theta\), but also on a large-scale variable
\(x=\varepsilon^2 \theta\). In this decomposition, \(\tilde b\in\mc
H_{\text{loc}}^s(\mb R^3\times \mb T^3)\) is periodic and average-free with
respect to the second variable \(\theta\), whereas \(\bar b\in\mc
H_{\text{loc}}^s(\mb R^3)\).

Using the formal variable separation above, the equation then writes under
the form
\begin{equation}
  \begin{split}
    \varepsilon\lambda^\varepsilon\tilde b &+ \lambda^\varepsilon\bar b
    - \nabla_x\wedge(U\wedge\tilde b) - \frac1\varepsilon\nabla_x\wedge(U\wedge\bar b) -\frac1{\varepsilon^2}\nabla_\theta\wedge(U\wedge\tilde b) \\
    &-\frac1{\varepsilon^3}\nabla_\theta\wedge(U\wedge\bar b)
    - \frac1{\varepsilon^3}\Delta_\theta\tilde b - \varepsilon\Delta_x\tilde b - \Delta_x\bar b - \frac2\varepsilon\sum_{i=1}^3\partial_{x_i\theta_i}\tilde b = 0.
  \end{split}
\end{equation}
We split this equation between mean part and fluctuating part:
\begin{equation} \label{sepmeanfluct}
  \begin{split}
    \Delta_\theta\tilde b &+ \varepsilon\nabla_\theta\wedge(U\wedge\tilde b) + \nabla_\theta\wedge(U\wedge\bar b) + \varepsilon^2(\nabla_x\cdot\nabla_\theta + \nabla_\theta\cdot\nabla_x)\tilde b + \\
    &\varepsilon^2\nabla_x\wedge(U\wedge\bar b) + \varepsilon^3\nabla_x\wedge\widetilde{(U\wedge\tilde b)} + \varepsilon^4\Delta_x\tilde b = \varepsilon^4\lambda^\varepsilon\tilde b,
  \end{split}
\end{equation}
\begin{equation}
  \label{eqn_mean}
  \nabla_x\wedge\overline{U\wedge\tilde b} + \Delta_x\bar b = \lambda^\varepsilon\bar b.
\end{equation}
And we rewrite the fluctuating part of the equation to separate terms with
\(\tilde b\) and \(\bar b\):
\begin{equation}
  \label{eqn_fluct}
  \begin{split}
    \Delta_\theta\tilde b &+\varepsilon\nabla_\theta\wedge(U\wedge\tilde b) + \varepsilon^2(\nabla_x\cdot\nabla_\theta + \nabla_\theta\cdot\nabla_x)\tilde b  +  \varepsilon^3\nabla_x\wedge\widetilde{(U\wedge\tilde{ b})} + \varepsilon^4\Delta_x\tilde b - \varepsilon^4\lambda^\varepsilon\tilde b\\
    &= -\nabla_\theta\wedge(U\wedge\bar b) - \varepsilon^2\nabla_x\wedge(U\wedge\bar b) 
  \end{split}
\end{equation}
Equations \eqref{sepmeanfluct} to \eqref{eqn_fluct} will be used
extensively throughout the note.

We can now give the theorems that we will prove in this note. The first one
is related to the induction equation \eqref{mhd}:
\begin{theoreme}
  Let \(s>0\) be a real number. There exists a dense subset \(\Omega\) of
  \(\ms P\) such that for all \(U\in\Omega\), there exists a real
  \(\varepsilon_0>0\), a vector \(\xi\in\mb R^3\) and a function
  \(\fonction{\lambda_U}[[0,\varepsilon_0[][\mb R_+^*]\) analytic such that
  for all \(\varepsilon\in[0,\varepsilon_0[\), there exists a solution
  \(b\) of \eqref{mhd} of the form:
  \begin{equation}
    b^\varepsilon(\theta,t) = e^{\lambda_U(\varepsilon) t + \varepsilon^2\iim\xi\cdot\theta}(\bar b' + \varepsilon\tilde b'(\theta)),
  \end{equation}
  where \(U\in\ms P\), \(\bar b'\in\mb R^3\) and \(\tilde b'\in\mc H_0^s(\mb T^3)\).
\end{theoreme}
Our second theorem concerns the full MHD system \eqref{fullmhd}.
For this next theorem, we introduce for any vector \(T\in(\mb R_+^*)^3\)
the Sobolev space \(\mc H^s_T = \mc H^s(\mb R/T_1\mb Z\times\mb R/T_2\mb
Z\times\mb R/T_3\mb Z)\)
\begin{theoreme}
  \label{lin_theoreme}
  Suppose \(s>3/2+1\). Then there exists a dense subset \(\Omega\) of \(\ms
  P\) such that for all \(U\in\Omega\), there exists a vector \(T =
  (T_1,T_2,T_3)\) with integer coordinates and a sequence
  \((\varepsilon_n)_{n\ge0}\) decreasing to \(0\) such that the solution
  \(\begin{pmatrix} U \\ 0\end{pmatrix}\) of the (MHD) system
  \eqref{inimhd} is nonlinearly unstable in \(\mc H_T^s\) in the following
  sense:
  
  There exists initial values
  \begin{equation}
    \begin{pmatrix} u_0 \\ b_0\end{pmatrix}\in\mc H_T^s,
  \end{equation}
  and \(C_0>0\) such that for all \(\delta>0\), the solution \(\begin{pmatrix} u_\delta \\ b_\delta \end{pmatrix}\) of \eqref{fullmhd} with initial value \(\delta\begin{pmatrix} u_0 \\ b_0\end{pmatrix}\) satisfies
  \begin{equation}
    \left\|\begin{pmatrix} u_\delta \\ b_\delta \end{pmatrix}(t_\delta)\right\|_{\mc H_T^s}\ge C_0
  \end{equation}
  for some time \(t_\delta\).
\end{theoreme}

This result extends substantially the article \cite{gerardvaret2006osi} by
D. G\'erard-Varet. Recast in our variables, this article considers the case
of large magnetic diffusion, that is \(R_m=\varepsilon\ll1\). Let us point
out that D.  G\'erard-Varet didn't use the physicists scaling for the
velocity field, thus making possible confusion to the fact that he chose
\og\(R_m = 1\)\fg{} and that it actually corresponds to the case where
\(R_m\) is vanishing with \(\varepsilon\). To get the usual Reynolds number
you have to do a time-space transformation \(t\mapsto t/\varepsilon^3\),
\(x\mapsto x/\varepsilon^2\) (and multiply the amplitude \(V\) by
\(1/\varepsilon\)), thus making the Reynolds number \(\varepsilon\) instead
of \(1\).

Briefly, it is shown in \cite{gerardvaret2006osi} that for all \(m\in\mb
N\), there are solutions of \eqref{mhd} (on a larger box depending on
\(\varepsilon\)) that go from amplitude \(\varepsilon^m\) initially to an
amplitude \(\eta = \eta(m)>0\) independent of \(\varepsilon\) (but
depending on \(m\)).  Note that, as \(\eta\) depends on \(m\), the result
in \cite{gerardvaret2006osi} does not yield Lyapunov instability at fixed
small \(\varepsilon\), contrary to ours.

The main idea to prove our two theorems is to construct an exact solution
\((\bar b, \tilde b)\) of equations \eqref{sepmeanfluct} to
\eqref{eqn_fluct}. We shall obtain it through a perturbative argument,
starting from the case \(\varepsilon = 0\).

\section{Case \texorpdfstring{$\varepsilon = 0$}{epsilon = 0}}
\subsection{An eigenvalue problem}
The fluctuating part in the case \(\varepsilon = 0\) is
\begin{equation}
  \Delta_\theta\tilde b = -\nabla_\theta\wedge(U\wedge\bar b).
\end{equation}
Since the Laplacian is invertible as a function from \(\mc H_0^{s+2}\) to
\(\mc H_0^s\), we can write
\begin{equation}
  \tilde b = -\Delta_\theta^{-1}\nabla_\theta\wedge(U\wedge\bar b) = \ms L(\theta)\bar b,
\end{equation}
where \(\fonction{\ms L(\theta)}{\mc H_0^s}{\mc H_0^{s+1}}\) (\(x\) can be
taken as a parameter since it is not involved in \(\ms L\)).

Now we can replace \(\tilde b\) by \(\ms L(\theta)\bar b\) in the mean part
of the equation \eqref{eqn_mean} to have an equation involving only \(\bar
b\):
\begin{equation}
  \nabla_x\wedge\overline{U\wedge\ms L(\theta)\bar b} + \Delta_x\bar b = \lambda^0\bar b.
\end{equation}
We introduce the matrix \(\alpha\) defined by
\begin{equation}
  \alpha(U)\bar b = \overline{U\wedge\ms L(\theta)\bar b} = -\int_{[0,1]^3}U\wedge\Delta_\theta^{-1}\nabla_\theta\wedge(U\wedge\bar b)\ud\theta.
\end{equation}
(The name of the matrix is the one related to the so called alpha effect
mentionned in the introduction)

Now the equation is
\begin{equation}
  \label{eq_eps0}
  \nabla_x\wedge(\alpha(U)\bar b) + \Delta_x\bar b = \lambda^0\bar b.
\end{equation}

We will find a solution of this last equation under the form \(\bar b(x) =
e^{\iim\xi\cdot x}\bar b'\). Using this stanza, the equation then has the
form (dropping the prime to lighten notations in the end of the
subsection):
\begin{equation}
  \label{eq_mean_eps0}
  \iim\xi\wedge(\alpha(U)\bar b) -|\xi|^2\bar b = \lambda^0\bar b.
\end{equation}

Thus we have come now to an eigenvalue problem with the matrix
\(\xi\wedge\alpha(U)\). The properties of the matrix \(\alpha\) are the
source of the alpha-effect mechanism; they have been already studied in
\cite{gerardvaret2006osi} and in \cite{bouya2011instability} with a more
general case. We give here a proposition which is only the particular case
of the similar study done in those papers:

\begin{proposition}
  The matrix \(\alpha(U)\) defined above is real and
  symmetric. Furthermore, there exists a dense subset \(\Omega\) of \(\ms
  P\) such that for all profile \(U\in\Omega\), there exists some
  \(\xi\in\mb R^3\) such that the matrix defined by the left hand side of
  \eqref{eq_mean_eps0} admits an eigenvalue of positive real part.
\end{proposition}
The reader should refer to those paper for a proof of this proposition.

\emph{From now on, we will assume that our base flow belongs to this
  subset \(\Omega\). This gives the existence of a vector \(\xi\in\mb R^3\)
  and an eigenvector \(b_0\) associated to an eigenvalue \(\lambda^0\) with
  positive real part such that
  \begin{equation}
    \bar b(x) = e^{\iim\xi\cdot x}b_0
  \end{equation}
  is a growing mode of \eqref{eq_eps0}}


\section{Existence of growing mode for \texorpdfstring{$\varepsilon>0$}{epsilon > 0}}
In this section we will search for a solution of the induction equation for
\(\varepsilon>0\), and we will search a solution under the form
\begin{equation}
  b^\varepsilon(\theta,t) = e^{\lambda^\varepsilon t + \varepsilon^2\iim\xi\cdot\theta}(\bar b + \varepsilon\tilde b(\theta)),
\end{equation}
based on the result of last section, where \(\xi\) is the one chosen for
\(\varepsilon = 0\). Thus the equations become:
\begin{equation}
  \label{eqn_eps}
  \begin{split}
    \Delta_\theta\tilde b &+\varepsilon\nabla_\theta\wedge(u\wedge\tilde b) + \iim\varepsilon^2(\xi\cdot\nabla_\theta + \nabla_\theta\cdot\xi)\tilde b  +  \varepsilon^3\iim\xi\wedge\widetilde{(u\wedge\tilde b)} - \varepsilon^4|\xi|^2\tilde b - \varepsilon^4\lambda^\varepsilon\tilde b\\
    &= -\nabla_\theta\wedge(u\wedge\bar b) - \varepsilon^2\iim\xi\wedge(u\wedge\bar b) ,
  \end{split}
\end{equation}
\begin{equation}
  \iim\xi\wedge\overline{u\wedge\tilde b} - |\xi|^2\bar b = \lambda^\varepsilon\bar b.
\end{equation}
    
The left part of \eqref{eqn_eps} may be seen as an operator from \(\mc
H_0^{s+2}\) to \(\mc H_0^s\), and for small enough \(\varepsilon\) and
\(\lambda^\varepsilon\) close enough to \(\lambda^0\) it is a perturbation
of the Laplacian. Thus this operator is invertible and we can write again
\begin{equation}
  \tilde b = \ms L_{\xi,\varepsilon,\lambda^\varepsilon}(\theta)\bar b,
\end{equation}
where \(\ms L_{\xi,\varepsilon,\lambda^\varepsilon}(\theta) =
-L_{\xi,\varepsilon,\lambda^\varepsilon}^{-1}\circ(\nabla_\theta\wedge
u\wedge\cdot+\varepsilon^2\iim\xi\wedge(u\wedge\cdot))\) is the operator
from \(\mc H_0^s\) to \(\mc H_0^{s+1}\) where
\begin{equation}
  L_{\xi,\varepsilon,\lambda^\varepsilon} = \Delta_\theta +\varepsilon\nabla_\theta\wedge(u\wedge\cdot) + \iim\varepsilon^2(\xi\cdot\nabla_\theta + \nabla_\theta\cdot\xi) +  \varepsilon^3\iim\xi\wedge\widetilde{(u\wedge\cdot)} - \varepsilon^4|\xi|^2 - \varepsilon^4\lambda^\varepsilon
\end{equation}
and the mean-part equation is
\begin{equation}
  \iim\xi\wedge(\alpha_{\xi,\varepsilon,\lambda^\varepsilon}(u)\bar b) - |\xi|^2\bar b = \lambda^\varepsilon\bar b,
\end{equation}
where the matrix \(\alpha\) depends here on some parameters:
\begin{equation}
  \begin{split}
    \alpha_{\xi,\varepsilon,\lambda^\varepsilon}(u)\bar b &= \overline{u\wedge\ms L_{\xi,\varepsilon,\lambda^\varepsilon}(\theta)\bar b} \\
    & = \int_{[0,1]^3}u\wedge L_{\xi,\varepsilon,\lambda^\varepsilon}^{-1}(\nabla_\theta\wedge u\wedge\bar b+\varepsilon^2\iim\xi\wedge(u\wedge\bar b))
  \end{split}
\end{equation}

This equation admits non-zero solution \(\bar b\) if and only if
\begin{equation}
  \det(\iim\xi\wedge\alpha_{\xi,\varepsilon,\lambda^\varepsilon}(u) - |\xi|^2\Id - \lambda^\varepsilon\Id) = 0,
\end{equation}
where
\begin{equation}
  \iim\xi\wedge\alpha_{\xi,\varepsilon,\lambda^\varepsilon}(u) = \iim\begin{pmatrix} 0 & -\xi_3 & \xi_2 \\ \xi_3 & 0 & -\xi_1 \\ -\xi_2 & \xi_1 & 0 \end{pmatrix}\alpha_{\xi,\varepsilon,\lambda^\varepsilon}(u) = A^\xi\alpha_{\xi,\varepsilon,\lambda^\varepsilon}(u).
\end{equation}

\(\xi\) being fixed by the former section, we define the following function
\begin{equation}
  f(\varepsilon,\mu) = \det(\iim\xi\wedge\alpha_{\xi,\varepsilon,\mu}(u) - |\xi|^2\Id - \mu\Id),
\end{equation}
defined in a neighborhood \(V\) of \((0,\lambda^0)\) in \(\mb R^2\).

\section{Implicit function theorem}
We know from former section that \(f(0,\lambda^0) = 0\). Now we want to
prove that for \(\varepsilon>0\) small enough, there exists
\(\lambda^\varepsilon\) such that \(f(\varepsilon,\lambda^\varepsilon) =
0\). Furthermore, \(f\) is analytic in \((\varepsilon,\mu)\) in \(V\).

We prove in this section that \(\left.\frac{\partial
    f}{\partial\mu}\right|_{\begin{subarray}{l} \varepsilon = 0\\ \mu =
    \lambda^0\end{subarray}}\neq 0\), so that we can apply the analytic
implicit functions theorem and conclude that there exists a neighborhood
\(W\) of \(0\) and a function \(\fonction{\varphi}{W}{\mb R}\) such that
\(\varphi(0) = \lambda^0\) and \(\forall\varepsilon\in
W,f(\varepsilon,\varphi(\varepsilon)) = 0\) (we can then restrict this
neighborhood to positive values of \(\varepsilon\)).

We compute now \(\left.\frac{\partial
    f}{\partial\mu}\right|_{\begin{subarray}{l} \varepsilon = 0\\ \mu =
    \lambda^0\end{subarray}}\): using usual computation rules on
composition of functions, we have that
\begin{equation}
  \left.\frac{\partial f}{\partial\mu}\right|_{\begin{subarray}{l} \varepsilon = 0\\ \mu = \lambda^0\end{subarray}} =
  \left.\ud\det\right|_{\iim\xi\wedge\alpha_{\xi,0,\lambda^0}(u)-|\xi|^2\Id - \lambda^0\Id}\circ\left.\ud(\iim\xi\wedge\alpha_{\xi,\varepsilon,\mu}(u)-|\xi|^2\Id - \mu\Id)\right|_{\begin{subarray}{l} \varepsilon = 0\\ \mu = \lambda^0\end{subarray}}(0,1).
\end{equation}
The computation of the second differential is straightforward and gives
\begin{equation}
  \left.\ud(\iim\xi\wedge\alpha_{\xi,\varepsilon,\mu}(u)-|\xi|^2\Id - \mu\Id)\right|_{\begin{subarray}{l} \varepsilon = 0\\ \mu = \lambda^0\end{subarray}}(0,1) = \left.\frac{\partial(\iim\xi\wedge\alpha_{\xi,\varepsilon,\mu}(u)-|\xi|^2\Id - \mu\Id)}{\partial\mu}\right|_{\begin{subarray}{l} \varepsilon = 0\\ \mu = \lambda^0\end{subarray}} = -\Id.
\end{equation}
And the computation of the first one gives
\begin{equation}
  \left.\ud\det\right|_{\iim\xi\wedge\alpha_{\xi,0,\lambda^0}(u)-|\xi|^2\Id - \lambda^0\Id}(-\Id) = -\Tr(\com(\iim\xi\wedge\alpha_{\xi,0,\lambda^0}(u)-|\xi|^2\Id - \lambda^0\Id)).
\end{equation}
We know from the second section the matrix \(A =
\iim\xi\wedge\alpha_{\xi,0,\lambda^0}(u)-|\xi|^2\Id - \lambda^0\Id\) is
diagonalizable and that his eigenvalues are
\(0,-\lambda^0,-2\lambda^0\). Thus, \(\com(A)\) is similar to the matrix
\begin{equation}
  \com\begin{pmatrix} 0 & 0 & 0 \\ 0 & -\lambda^0 & 0 \\ 0 & 0 & -2\lambda^0\end{pmatrix} = \begin{pmatrix} 2(\lambda^0)^2 & 0 & 0 \\ 0 & 0 & 0 \\ 0 & 0 & 0\end{pmatrix},
\end{equation}
which is of trace \(2(\lambda^0)^2\neq 0 \).

Thus
\begin{equation}
  \left.\frac{\partial f}{\partial\mu}\right|_{\begin{subarray}{l} \varepsilon = 0\\ \mu = \lambda^0\end{subarray}} = 2(\lambda^0)^2\neq 0 .
\end{equation}
Since \(\varphi\) is continuous on \(0\), and \(\varphi(0) = \lambda^0\)
has positive real value, this remains true for \(\varepsilon>0\)
sufficiently small, which concludes the proof of the existence of a growing
mode for a small \(\varepsilon>0\).

\section{Nonlinear instability}
We will now prove nonlinear instability of the system \eqref{mhd}. For that
we will adapt a method developped by S. Friedlander, W. Strauss and
M. Vishik in \cite{friedlander1997nonlinear} (itself adapted from Y. Guo
and W. Strauss in \cite{guo1995nonlinear}) to prove nonlinear instability
given linear instability. The condition \(s>3/2+1\) in the main theorem is
used only in this section to obtain an inequality on the nonlinear terms
similar to that in the article \cite{friedlander1997nonlinear}.  Actually,
using the parabolic feature of the MHD system, it would be possible to
lower our regularity requirements (see for instance Friedlander et al
\cite{friedlander2006nonlinear} in the context of the Navier-Stokes
equations). We stick here to regular data for simplicity of exposure.

Let \(U\in\ms P\) belonging to the dense subset of linearly unstable flows
described in the previous section.

We introduce Leray operator \(\mb P\) to get rid of \(p\). We rewrite
system \eqref{fullmhd} as
\begin{equation}
  \label{MHDPert2}
  \partial_t\ms U + L_s\ms U = Q(\ms U,\ms U),
\end{equation}
where
\begin{equation}
  \ms U = \begin{pmatrix} u' \\ b \end{pmatrix},
\end{equation}
\begin{equation}
  L_s\ms U = \varepsilon^{-3}\begin{pmatrix} \mb P(u\cdot\nabla u' + u'\cdot\nabla u) - \Delta u' \\
    -\nabla\wedge(u\wedge b) - \Delta b\end{pmatrix}
\end{equation}
is the linearized operator around \(\begin{pmatrix} u\\0\end{pmatrix}\),
and
\begin{equation}
  Q(\ms U,\ms U) = \varepsilon^{-3}\begin{pmatrix} \mb P((\nabla\wedge b)\wedge b - u'\cdot\nabla u') \\ \nabla\wedge(u'\wedge b)\end{pmatrix}
\end{equation}
is the nonlinearity.

We shall consider this system in \(\mc H_T^s\), \(T\) to be specified
later. We also denote \(L^2_T = \mc H_T^0\) the \(L^2\) space on the
torus. Since we will work in these spaces until the end of the proof, we
will drop the \(T\) in the notation. We have the following a priori
estimates for the nonlinear term:

Every term of every component of \(Q(\ms U,\ms U)\) can be written as a
product of one component of \(\ms U\) and one of \(\nabla\ms U\) (composed
by \(\mb P\)). Thus, we have the first estimate
\begin{equation}
  \|Q(\ms U,\ms U)\|_{L^2}\le C\varepsilon^{-3} \||\ms U||\nabla\ms U|\|_{L^2} \le C\varepsilon^{-3} \|\ms U\|_{L^2}\|\nabla\ms U\|_{L^\infty}.
\end{equation}
For \(r = \frac12(s+n/2+1) > n/2+1\) (here \(n=3\)), we have that
\begin{equation}
  \|\nabla\ms U\|_{L^\infty}\le \|\ms U\|_{\mc H^r},
\end{equation}
and
\begin{equation}
  \|\ms U\|_{\mc H^r} \le C \|\ms U\|_{L^2}^\eta\|\ms U\|_{\mc H^s}^{1-\eta}
\end{equation}
for \(\eta\) such that \(r = (1-\eta)s\) (i.e. \(\eta = \frac12 -
\frac{n+2}{4s} = \frac12 - \frac5{4s}\)).

Thus, we have:
\begin{equation}
  \|Q(\ms U,\ms U)\|_{L^2}\le C\varepsilon^{-3} \|\ms U\|_{L^2}^{1+\eta}\|\ms U\|_{\mc H^s}^{1-\eta}
\end{equation}

We shall prove the following theorem, which is a mere reformulation of the
second theorem:
\begin{theoreme}
  Suppose \(s>3/2+1\). Then there exists a vector \(T = (T_1,T_2,T_3)\)
  with integer coordinates and a sequence \((\varepsilon_n)_{n\ge0}\)
  decreasing to \(0\) such that the system \eqref{MHDPert2} is nonlinearly
  unstable in \(\mc H_T^s\) in the following sense:

  There exists a growing mode \(\ms U_0=\begin{pmatrix} u_0 \\
    b_0\end{pmatrix}\) with \(\|\ms U_0\|_{\mc H_T^s} = 1\) of the
  linearized system and a constant \(C_0>0\) depending only on \(\ms
  U_0,s,n,\rho\) such that for all \(\delta>0\), the solution of the full
  magnetohydrodynamics system with initial value \(\delta\ms U_0\)
  satisfies \(\|\ms U(t_\delta)\|_{\mc H_T^s}\ge C_0\) for some time
  \(t_\delta\).
\end{theoreme}
\begin{remarque}
  This theorem only proves that the solution reaches high values in finite
  time no matter how small the initial condition is. In particular it doesn't say
  wether it explodes in finite time or is defined for all time, since both
  situations are possible in the result of the theorem. Actually, for small
  enough initial data (that is small enough \(\delta\)), it is classical
  that the solution of the system doesn't explode in finite time. More
  details in \cite{temam2001navier} (Theorem 3.8)
\end{remarque}

Assume in contrary that the system is nonlinearly stable while spectrally
unstable, that is that for any growing mode of initial value \(\ms
U_0=\begin{pmatrix} u_0 \\ b_0\end{pmatrix}\) with \(\|\ms U_0\|_{\mc
  H_T^s} = 1\) in any box \(T\) of the system (i.e. an eigenvector of
\(L_s\) associated to an eigenvalue with positive real part) and for any
\(\varepsilon>0\), there exists a \(\delta>0\) such that for all \(t>0\),
the solution \(\ms U(t)\) of the system with initial values \(\delta\ms
U_0\) satisfies \(\|\ms U(t)\|_{\mc H_T^s}\le\varepsilon\).

\subsection{Proof of the nonlinear instability}
Let \(b_L(\theta) =
e^{\iim\varepsilon^2\xi\cdot\theta}b^\varepsilon(\theta)\) be the unstable
eigenmode given by Theorem \(2\). We fix \(\varepsilon<\varepsilon_0\) in
such a way that
\begin{equation}
  T = \left(\frac{2\pi}{\varepsilon^2|\xi_1|},\frac{2\pi}{\varepsilon^2|\xi_2|},\frac{2\pi}{\varepsilon^2|\xi_3|}\right)
\end{equation}
belongs to \(\mb N^3\). This is possible, as it is clear from the proof of
Theorem \ref{lin_theoreme} that \(\xi\) can be chosen in \((2\pi\mb
Q_+^*)^3\) (take \(\varepsilon\) to be the inverse of a large integer,
which can be chosen as big as we want). In particular, \(\tilde{\ms U_0}
= \begin{pmatrix} 0\\b_L\end{pmatrix}\) is an eigenvector of \(L_s\)
associated to an eigenvalue with positive real part in \(\mc H_T^s\).

Denote by \(\rho\) the maximum real part of the spectrum of \(L_s\):
\begin{equation}
  \rho = \max\left\{\Re\lambda,\lambda\in\spectr(L_s)\right\}.
\end{equation}
By the previous remark on \(\tilde{\ms U_0}\), we know that
\(\rho>0\). Moreover, since the spectrum of \(L_s\) is only made of
eigenvalues, there exists an eigenvector \(\ms U_0 = \begin{pmatrix}
  u_0\\b_0\end{pmatrix}\) with eigenvalue \(\lambda\) satisfying exactly
\(\Re\lambda = \rho\).

Moreover, by standard properties of the spectral radius of \(e^{L_s}\), we
know that
\begin{equation}
  e^\rho = \lim_{t\rightarrow+\infty}\|e^{tL_s}\|^{1/t}.
\end{equation}
Thus, for all \(\eta>0\), there exists a \(C_\eta\) such that for all
\(t\ge0\),
\begin{equation}
  \|e^{tL_s}\|\le C_\eta e^{\rho(1+\frac\eta2)t}
\end{equation}

Then, the solution \(\ms U\) of \eqref{MHDPert2} with initial data
\(\delta\ms U_0\) can be written
\begin{equation}
  \ms U = \dot{\ms U}+\delta e^{\lambda t}\ms U_0,
\end{equation}
where \(\dot{\ms U}\) satisfies
\begin{equation}
  \left\{\begin{array}{l}
      \partial_t\dot{\ms U} = L_s\dot{\ms U} + Q(\ms U,\ms U)\\
      \dot{\ms U}_{|t=0}=0
    \end{array}\right.
\end{equation}

Using Duhamel's formula, we can write the solution under the form:
\begin{equation}
  \dot{\ms U}(t) = \int_0^te^{(t-s)L_s}Q(\ms U,\ms U)\ud s.
\end{equation}
	
Define
\begin{equation}
  T_\delta = \sup\left\{t>0,\forall s\le t,\|\dot{\ms U}(s)\|_{L^2}\le\frac12\delta e^{\rho s}\|\ms U_0\|_{L^2}\right\},
\end{equation}
where \(\lambda\) is the coefficient in the growing linear mode. Then,
since we supposed that \(\ms U = \ms U_l + \dot{\ms U}\) remains bounded,
and \(\|\ms U_l\|\sim \delta e^{\rho t}\|\ms U_0\|\), we have that
\(T_\delta < \infty\), and for all \(t<T_\delta\),
\begin{equation}
  \|\ms U(t)\|_{L^2}\le \frac32\delta e^{\rho t}\|\ms U_0\|_{L^2}.
\end{equation}

Using a priori estimates stated earlier, we have for \(\eta = \frac12 -
\frac{n+2}{4s} = \frac12 - \frac{5}{4s}\):
\begin{equation}
  \|\dot{\ms U}(t)\|_{L^2} \le C_\eta\varepsilon^{-3} \int_0^te^{\rho(1+\frac\eta2)(t-s)}\|\ms U\|_{L^2}^{1+\eta}\|\ms U\|_{\mc H^s}^{1-\eta}\ud s.
\end{equation}
\begin{equation}
  \begin{split}
    \|\dot{\ms U}(t)\|_{L^2} &\le \delta^{1+\eta}C_\eta\varepsilon^{-3} \int_0^te^{(t-s)(1+\frac\eta2)\rho}e^{(1+\eta)\rho s}\varepsilon^{1-\eta}\ud s \\
    &\le C_\eta\delta^{1+\eta}e^{t(1+\frac\eta2)\rho}\varepsilon^{-2-\eta}\frac2{\eta\rho}(e^{t\frac\eta2\rho}-1) \\
    &\le C'\delta^{1+\eta}e^{t(1+\eta)\rho}\varepsilon^{-2-\eta}\\
  \end{split}
\end{equation}
Thus, getting back to \(\ms U\):
\begin{equation}
  \|\ms U\|_{L^2}\ge\delta\|\ms U_0\|_{L^2}e^{\rho t} - C'\delta^{1+\eta}\varepsilon^{-2-\eta}e^{t(1+\eta)\rho}.
\end{equation}
By definition of \(T_\delta\), we have:
\begin{equation}
  \|\dot{\ms U}(T_\delta)\|_{L^2} = \frac12\delta e^{\rho T_\delta}\|\ms U_0\|_{L^2}\le C'\delta^{1+\eta}\varepsilon^{-2-\eta}e^{T_\delta(1+\eta)\rho},
\end{equation}
that is
\begin{equation}
  \frac12\|\ms U_0\|_{L^2}\le C'\delta^\eta\varepsilon^{-2-\eta}e^{T_\delta\eta\rho},
\end{equation}
\begin{equation}
  T_\delta\ge\ln\left(\frac{\|\ms U_0\|_{L^2}}{2C'\delta^\eta\varepsilon^{-2-\eta}}\right)\frac1{\eta\rho}.
\end{equation}
Denote by \(t_\delta\) the quantity
\begin{equation}
  t_\delta = \ln\left(\frac{\|\ms U_0\|_{L^2}}{2C'\delta^\eta\varepsilon^{-2-\eta}(1+\eta)}\right)\frac1{\eta\rho}.
\end{equation}
Thus \(t_\delta \le T_\delta\), and at this time \(t_\delta\)
\begin{equation}
  \|\ms U\|_{L^2} \ge \delta e^{\rho t_\delta}\left(\|\ms U_0\|_{L^2} - \frac{\|\ms U_0\|_{L^2}}{1+\eta}\right) = \frac\eta{1+\eta}\|\ms U_0\|_{L^2}\left(\frac{\|\ms U_0\|_{L^2}\varepsilon^{2+\eta}}{C'(1+\eta)}\right)^{1/\eta}.
\end{equation}
Thus \(\|\ms U\|_{L^2}\) (and \(\|\ms U\|_{\mc H^s}\)) can be bounded below
by a value independent of \(\delta\), which contradicts the assumption, and
proves the theorem.

\bibliographystyle{abbrv}
\bibliography{biblio_perso}
\end{document}